Volodymyr Denysiuk, Dr.Phys.-Math.Sci., Professor
(National Aviation University, Ukraine)

# METHOD OF IMPROVEMENT OF CONVERGENCE FOURIER SERIES AND INTERPOLIATION POLYNOMIALS IN ORTHOGONAL FUNCTIONS

**Abstract**

There is proposed a method for improving the convergence of Fourier series by function systems, orthogonal at the segment, the application of which allows for smooth functions to receive uniformly convergent series. There is also proposed the method of phantom nodes improving the convergence of interpolation polynomials on systems of orthogonal functions, the application of which in many cases can significantly reduce the interpolation errors of these polynomials. The results of calculations are given at test cases using the proposed methods for trigonometric Fourier series; these calculations illustrate the high efficiency of these methods. Undoubtedly, the proposed method of phantom knots requires further theoretical studies.

**Introduction.**

In many science and technology tasks arises a problem of representation of a function given on a segment $[\alpha,\beta]$, uniformly convergent Fourier series for some complete system of functions, orthogonal at a section $[a,b]$ with some weight. Since the given segment of the investigated function $[\alpha,\beta]$ does not coincide with the segment $[a,b]$, then the linear replacement of the variables segment $[\alpha,\beta]$ is reflected on the segment of the orthogonality $[a,b]$ of the selected system of functions. This case is well investigated in the theory of approximation.

However, in many cases the obtained Fourier series in this case coincides uniformly only in the interval $[a_1, b_1] \in (a,b)$ [1]; this situation is explained by the properties of the selected system of orthogonal functions with a given weight function.

For example, if the function $f(x)$, $x \in [\alpha,\beta]$ has a continuous derivative of the third order, then its Fourier series in polynomials of Chebyshev of the second kind coincides uniformly on any interval $[-1+h, 1-h]$, $0 < h < 1$ [2]. Another example is the trigonometric functions system, since in the expansion of a smooth function $f(x)$, $x \in [\alpha,\beta]$ in a Fourier trigonometric series, this function must be continued periodically for the entire numerical axis; it is clear that with such a continuation, in the general case In the periodic function and its derivatives there are even ruptures of the first kind of jump type in points $a,b$ and the uniform convergence of a row is provided at any segment $[a_1, b_1] \in (a,b)$.

**Formulation of the problem.**

Development and research of the method of improving the convergence of Fourier series in complete systems of orthogonal functions; development and research of the discrete version of the proposed method for improving the convergence of Fourier series.

**Main part.**

1. A method of improving the convergence of Fourier series in orthogonal systems of functions.

Let's consider a general method of obtaining uniformly convergent Fourier series in orthogonal systems at the segment $[a,b]$ functions for the above cases. This method is that a function $f(x)$ specified at a segment $[\alpha,\beta]$, the linear replacement of the variables is reflected at the segment $[a+p, b-q]$, $0 < p, q < (b-a)/2$; on segments $[a, a+p]$ and $(b-q, b]$ function $f(x)$ is completed by some functions $\varphi_1(x)$, $x \in [a, a+p)$ and $\varphi_2(x)$, $x \in (b-q, b]$, which are chosen in such a way to preserve that "good" property of the original function $f(x)$, which ensures uniform convergence of the Fourier series. In other words, in the Fourier series on the segment the function is now decomposed, which looks like

$$F(x) = \begin{cases} \varphi_1(x), & x \in [a, a+p); \\ f(x), & x \in [a+p, b-q]; \\ \varphi_2(x), & x \in (b-q, b]. \end{cases}$$

In many cases the parameters can be set equal, i.e. $p = q = h$; then uniform convergence to the function $f(x)$ is provided at the segment $[a+h, b-h]$; the violation of uniform convergence at the edges of the segment $[a, b]$ will only affect functions $\varphi_1(x)$, $x \in [a, a+h)$ and $\varphi_2(x)$, $x \in (b-h, b]$.

Further functions $\varphi_1(x)$ and $\varphi_2(x)$ will be called phantom functions, method of extending of function $f(x)$ functions $\varphi_1(x)$ and $\varphi_2(x)$ - the method of phantom functions.

Often, the property of a function $f(x)$ that ensures uniform convergence of the Fourier series in certain systems of orthogonal functions is that this function must be continuous and have continuous derivatives of certain orders. In such cases in the role of phantom functions, it is expedient to choose interpolation polynomials of Hermite, thanks to which it is possible to provide the necessary smoothness of function $F(x)$ on the segment $[a, b]$.

Let's consider the using of the method of phantom functions more detail, at the same time, without losing the universality, we confine ourselves to certainty only by considering well-studied trigonometric Fourier series.

Let $f(x)$ - a piecewise differentiated function given on a finite segment $[0, 2\pi]$; as is known, such functions may have a finite number of discontinuities of both the function itself and its derivatives. Assume that the functions $f, f', ..., f^{(k-1)}$ are continuous at $(0, 2\pi)$, excluding points $\xi_\mu, (\mu = 1, 2, ..., m)$, where all these functions have jumps, respectively equal $\delta_\mu^{(0)}, \delta_\mu^{(1)}, ..., \delta_\mu^{(k-1)}$. To these points you need to add a point $\xi_0 = 0$, if difference

$$\delta_0^{(0)} = f(0+0) - f(2\pi - 0);$$
$$\delta_0^{(1)} = f'(0+0) - f'(2\pi - 0);$$
$$\ldots$$
$$\delta_0^{(k-1)} = f^{(k-1)}(0+0) - f^{(k-1)}(2\pi - 0);$$

different from zero, since there are gaps in the periodic extension of the function at this point. Suppose also that there is a derivative everywhere (except points $\xi_\mu, (\mu = 1, 2, ..., m)$, there is a derivative $f^{(k)}$ that at $[-\pi, \pi]$ is a function of limited variation. Let's denote

$$A_n^{(i)} = -\frac{1}{\pi} \sum_{\mu=1}^{m} \delta_\mu^{(i)} \sin n\xi_\mu \; ; \quad B_n^{(i)} = \frac{1}{\pi} \sum_{\mu=0}^{m} \delta_\mu^{(i)} \cos n\xi_\mu \; .$$

Then for the Fourier coefficients $a_n$ and $b_n$, that are calculated by the formulas

$$a_n = \frac{1}{\pi} \int_0^{2\pi} f(x) \cos nx \, dx \; ; \quad b_n = \frac{1}{\pi} \int_0^{2\pi} f(x) \sin nx \, dx, \quad n = 1, 2, \ldots \; .$$

Formulas have a place (respectively, for $k$ odd and $k$ pair):

$$a_n = \frac{A_n^{(0)}}{n} - \frac{B_n^{(1)}}{n^2} - \frac{A_n^{(2)}}{n^3} + \frac{B_n^{(3)}}{n^4} + \ldots \begin{cases} + (-1)^{\frac{k-1}{2}} \frac{A_n^{(k-1)}}{n^k} + O\left(\frac{1}{n^{k+1}}\right), \\ + (-1)^{\frac{k}{2}} \frac{B_n^{(k-1)}}{n^k} + O\left(\frac{1}{n^{k+1}}\right). \end{cases}$$

$$b_n = \frac{B_n^{(0)}}{n} + \frac{A_n^{(1)}}{n^2} - \frac{B_n^{(2)}}{n^3} - \frac{A_n^{(3)}}{n^4} + \ldots \begin{cases} + (-1)^{\frac{k-1}{2}} \frac{B_n^{(k-1)}}{n^k} + O\left(\frac{1}{n^{k+1}}\right), \\ + (-1)^{\frac{k}{2}} \frac{A_n^{(k-1)}}{n^k} + O\left(\frac{1}{n^{k+1}}\right). \end{cases}$$

These formulas are easy to obtain, integrating expressions for Fourier coefficients in parts $k-1$ times and giving the obtained integral of the Styltes through the Riemann integrals.

It follows from the above formulas, that the presence of gaps both in the function itself and in its derivatives, reduces the order of decreasing of the Fourier coefficients and, accordingly, worsens the approximative properties of the Fourier series. In particular, it is clear that the Fourier series coincides evenly with the function $f(x)$, only if $A_n^{(0)}$ and $B_n^{(0)}$ equal to zero.

However, as we have already mentioned even if the function $f(x)$, given at the segment $[0, 2\pi]$, is continuous and has continuous derivatives $k-1$ in order inclusive at this interval, then with a periodic extension with a period $2\pi$ for the entire numerical axis in the general case there are gaps of the first type of jump type both in the function itself and its derivatives. It is clear that in this case, the Fourier series can not match evenly on this segment.

To obtain a uniformly convergent Fourier series of function $f(x)$ on a segment $[0, 2\pi]$, we apply the above method of phantom functions.

The linear replacement of a variable will display the function $f(x)$, $x \in [0, 2\pi]$, on the segment $x \in [\alpha/2, 2\pi - \alpha/2]$, $0 < \alpha < 2\pi$; at intervals $[0, \alpha/2]$ and $(2\pi - \alpha/2, 2\pi)$ the function is defined so that, with a periodic extension, the continuity of the function itself and its derivatives to a certain order is preserved. It's easy to achieve using Hermite interpolation polynomials. Thus, from the consideration of the function $f(x)$, $x \in [0, 2\pi]$ we switched to the consideration of the function

$$F(x) = \begin{cases} \varphi_1(x), & x \in [0, \alpha/2); \\ f(x), & x \in [\alpha/2, 2\pi - \alpha/2]; \\ \varphi_2(x), & x \in (2\pi - \alpha/2, 2\pi]. \end{cases}$$

which, after periodic continuation with a period $2\pi$ for the entire numerical axis, is continuous and has continuous derivatives of certain orders.

It is advisable to make the following remark. Since the periodic function with the period $2\pi$ can be considered at any length $2\pi$, after replacing the variable let's proceed to the consideration of the function

$$F(t) = \begin{cases} f(t), & t \in [0, 2\pi - \alpha]; \\ \varphi(t), & t \in (2\pi - \alpha, 2\pi]. \end{cases}$$

The expediency of such a transition is due to the fact that now, instead of two phantom functions $\varphi_1(x)$, $x \in [0, \alpha/2)$ and $\varphi_2(x)$, $x \in (2\pi - \alpha/2, 2\pi]$, we are dealing with only one phantom function $\varphi(t), t \in (2\pi - \alpha, 2\pi]$.

To illustrate the foregoing, consider an example.

Example. On the segment $[0, 2\pi]$ we will consider the continuous function

$$f(x) = x - \pi, \quad x \in [0, 2\pi]$$

The Fourier series of this function has such form

$$S(x) = -2 \sum_{n=1}^{\infty} \frac{\sin nx}{n}$$

Here, in pure form, there is a role of the function breaks that are formed at points $0$ and $2\pi$ with the periodic extension of the function $f(x)$; in particular, the presence of these discontinuities leads to the fact that the Fourier coefficients of this function are decreasing order $O(1/n)$, and the row does not coincide evenly.

To obtain a uniformly convergent series, we apply the proposed method of phantom functions. Let's set the parameter $\alpha$, ($0 < \alpha < 2\pi$) and we will display the function at the segment $[0, 2\pi - \alpha]$. It is easy to do this by linearly replacing a variable $x$ with a variable $t$ in this way

$$t = \frac{2\pi - \alpha}{2\pi} x.$$

We construct the function $\lambda(t)$ in this way

$$\lambda(t) = f(0) + \frac{f(2\pi) - f(0)}{\alpha} (2\pi - t), \quad t \in (2\pi - \alpha, 2\pi)$$

Consider the function

$$\varphi(t,\alpha) = \begin{cases} f(t), & t \in [0, 2\pi - \alpha] \\ \lambda(t), & t \in (2\pi - \alpha, 2\pi) \end{cases}$$

We will continue to function $\varphi(t,\alpha)$, $t \in [0, 2\pi)$ with a period $2\pi$ for the entire numerical axis. It is easy to see that the resulting periodic function is continuous and satisfies the Lipschitz condition of uniform convergence on any segment of a numerical axis (as a piecewise differentiating function); so its Fourier series coincides evenly, and at the segment $[0, 2\pi - \alpha]$ coincides evenly to the function $f(t)$, because $\varphi(t,\alpha) \equiv f(t)$, $t \in [0, 2\pi - \alpha]$. The Fourier series has functions that have the form

$$S(t,\alpha) = -2\sum_{n=1}^{\infty} \left[ \frac{\sin nt}{n} \cdot \frac{2\pi}{(2\pi - \alpha)} \frac{\sin(n\frac{\alpha}{2})}{n\frac{\alpha}{2}} \right].$$

It is easy to see that the introduction of the phantom function $\lambda(t)$ in order to eliminate gaps has led to the fact that the Fourier coefficients of functions $f(x)$ are multiplied by a linear factor

$$\sigma(n,\alpha) = \frac{2\pi}{(2\pi - \alpha)} \frac{\sin(n\frac{\alpha}{2})}{n\frac{\alpha}{2}}$$

It is clear that the introduction of such a factor in the Fourier series leads to the fact that the Fourier coefficients are decreasing $O(n^{-2})$; consequently, at the basis of Weierstrass, a number coincides uniformly throughout the segment $[0, 2\pi]$. We note that in more detail the linear method of summation of Fourier series with a multiplier $\sigma(k,\alpha)$ is considered in [3].

Thus, the elimination of the function breaks by the proposed method allowed us to obtain a uniformly convergent Fourier series, which in the segment $[0, 2\pi - \alpha]$ coincides with the function $f(t)$.

In case of a function $f(x) \in C^{k-1}_{[0,2\pi]}$, the obtained result can be greatly enhanced by ensuring the continuity of derivatives of higher orders. Construct a function based on conditions

$$\lambda(2\pi - \alpha) = f(2\pi); \lambda(2\pi) = f(0);$$
$$\lambda'(2\pi - \alpha) = f'(2\pi); \lambda'(2\pi) = f'(0);$$
$$\ldots\ldots\ldots\ldots\ldots\ldots\ldots\ldots\ldots\ldots$$
$$\lambda^{(k-1)}(2\pi - \alpha) = f^{(k-1)}(2\pi); \lambda^{(k-1)}(2\pi) = f^{(k-1)}(0).$$

As we have already said, it is easy to construct such a function using the Hermite interpolation polynomial. Considering now the function

$$\phi(t,\alpha) = \begin{cases} f(t), & t \in [0, 2\pi - \alpha] \\ \lambda(t), & t \in (2\pi - \alpha, 2\pi) \end{cases} \quad (1)$$

and continuing it periodically with a period on the entire numerical axis; It is easy to see that the function itself and its derivatives to $k-1$ are inclusive continuous on the whole number axis; consequently, the breaks can have only derivative $k$ order. According to formulas (1), the Fourier coefficients $a_n$, $a_n$ in this case have a decreasing order $O(n^{-(k+1)})$.

Note that the proposed method of phantom functions can be considered a generalization of the well-known O.S.Maliyev's method.

2. **The method of phantom nodes** of improving the convergence of interpolation trigonometric polynomials.

Consider the function

$$f(t) = t + 1, \quad t \in [0, 2\pi].$$

Let's analyze the function $f(t)$ with the step $h_k = 2\pi(i-1)/N$ putting, for example, $N = 9$. Calculating the value of the function, which in our case will be equal 1, 2, 3, 4, 5, 6, 7, 8, 9, we will construct a trigonometric interpolation polynomial. Here it is advisable to make a such remark. Generally speaking, an interpolation trigonometric polynomial $T_N(t)$ interpolates a function $f(t)$ in $N+1$ points given on a segment $2\pi$. The fact is that at the point $2\pi$ the value of this polynomial is definite due to periodicity $T_N(2\pi) = T_N(0)$. Therefore, the polynomial interpolation $T_N(t)$ will be considered only at the interval of interpolation $2\pi - h$, pointing this section every time. It's easy to calculate that the normalized interpolation error $\varepsilon(t)$, which is calculated by the formula $\varepsilon(t) = (f(t) - T_8(t))/|f_N - f_1|$, equals 165.

Let's analyze the result. It is clear that interpolation requires a trigonometric polynomial periodic extension of the function $f(t)$ to the entire numerical axis with the period $2\pi$; consequently, as before, at the points $2\pi k$ ($k = 0, \pm 1, \pm 2, \ldots$) there are gaps in the functions of the jump type. The presence of such discontinuities, in turn, leads to relatively large interpolation errors of the function $f(t)$.

So we come to the conclusion, that a trigonometric interpolation polynomial has the same defects as the Fourier series considered above; and in this case we will propose a method of improving convergence, which we will call the method of phantom knots.

This method is as follows. Add to the sequence of interpolation nodes on the right side a pair of nodes that we will call phantom.

It is clear that the addition $2k$ ($k = 1, 2, \ldots$) of phantom nodes increases the number of interpolation nodes per segment $[0, 2\pi]$ and reduces the step of an interpolation grid, which now becomes equal $h_k = 2\pi(i-1)/(N+2k)$. Since the number of values of interpolated function does not change then reducing the step of the interpolation grid leads to a decrease in the interval of interpolation, which becomes equal $Nh_k$. At the segment $[2\pi - Nh_k, 2\pi]$ we will build a function $\varphi(t)$ that will satisfy the conditions

$$\varphi(t) = \begin{cases} f_N, & t = Nh_k; \\ f_1, & t = 2\pi. \end{cases}$$

When calculating the value of this function in phantom nodes, we find values in these nodes. Since the addition of phantom nodes leads to an increase in the total number of interpolation points, which, in turn, leads to an increase in the order of the interpolation polynomial you can expect a significant reduction in the interpolation error on the segment $[0, 2\pi - Nh_k]$.

Note that the even number of points is selected only to ensure that the total number of interpolation nodes is odd, since the trigonometric interpolation of the odd number of points is more convenient. In our opinion, it is convenient to add a small number of phantom nodes, that is, to lay $k = 1, 2$.

When constructing a function it is possible to require that its derivatives of a certain order also take certain values at points $Nh_k$ and $2\pi$; to find these values, you can use the divided differences of the interpolated function or exact values of derivatives of an interpolated function $f(t)$, if it known.

Note that values in phantom nodes can be selected, based on other considerations, thus avoiding the need to build a function $\varphi(t)$.

It is easy to see that the method of phantom nodes is actually a discrete version of the method of phantom functions for improving the convergence of Fourier series, which we considered earlier. We will illustrate the effectiveness of the proposed method of phantom nodes, continuing to consider the function

$$f(t) = t + 1, \quad t \in [0, 2\pi].$$

To 9 nodes of interpolation add 2 phantom nodes; The step of the interpolation grid now is $h = 2\pi/11$. The values in phantom nodes will be calculated using the values of the interpolated function and its first divided divisions; consequently, the values of the interpolated sequence have the form

1, 2, 3, 4, 5, 6, 7, 8, 9, 7.053, 2.947

Here is a graph of an interpolation polynomial constructed by 11 points (Fig. 1) and the graph of the normalized interpolation error at the interval of interpolation $[0, 2\pi - 3h]$ (Fig. 2).

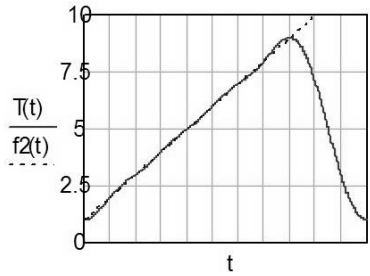

Fig. 1. Graph of polynomial $T(t)$ and function $f(t)$ on a segment $[0, 2\pi]$.

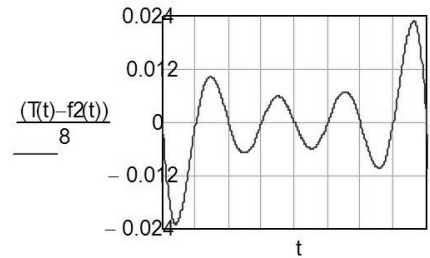

Fig. 2. Graph of the normalized error of interpolation with two phantom nodes on the interval of interpolation $[0, 2\pi - 3h]$; the error has decreased in $6.875$ times.

Add 4 phantom points now; the step of an interpolation grid now is $2\pi/13$. As before, the values in phantom nodes will be calculated using the values of the interpolated function and its first derivative; consequently, the value of the interpolated sequence has the form

1, 2, 3, 4, 5, 6, 7, 8, 9, 8.4, 6.3, 3.7, 1.6

Here is a graph of an interpolation polynomial constructed by 13 points (Fig. 3), and the graph of the normalized interpolation error at the interval of interpolation $[0, 2\pi - 5h]$ (Fig. 4).

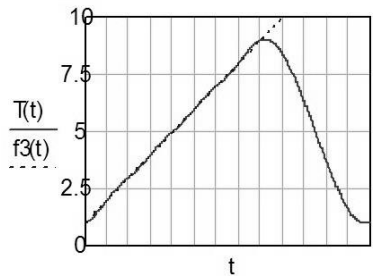

Fig. 3. Graph of polynomial $T(t)$ and function $f(t)$ on a segment $[0, 2\pi]$; interval of interpolation $[0, 2\pi - 5h]$.

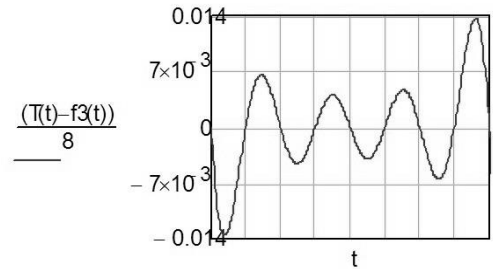

Fig. 4. Graph of the normalized error of interpolation with four phantom nodes on the interval of interpolation $[0, 2\pi - 5h]$; the error has decreased in $11.786$ times.

The results of further calculations are given in tables 1-6.

*Error tables for interpolation of function*
$$y = t+1, \; x \in [0, 2\pi)$$

Table 1. The ratio of interpolation errors without phantom nodes to the errors of interpolation with 2 phantom nodes

| Number of Interpolation points | Continuity of the function | Continuity of the function and the first derivative | Continuity of the function and the second derivative | Values are selected |
|---|---|---|---|---|
| 5 | 2.5 | 7.6 | 7.6 | 35.2 |
| 9 | 3 | 8,7 | 18.3 | 36.7 |
| 13 | 3.2 | 9.4 | 9.4 | 35.5 |

Table 2. The ratio of interpolation errors without phantom nodes to the errors of interpolation with 4 phantom nodes

| Number of Interpolation points | Continuity of the function | Continuity of the function and the first derivative | Continuity of the function and the second derivative | Values are selected |
|---|---|---|---|---|
| 5 | 3.4 | 15.8 | 41.3 | 208 |
| 9 | 4.5 | 20.6 | 61.1 | 1153.8 |
| 13 | 4.8 | 22.9 | 71.1 | 1159.4 |

Here are arrays of values in the interpolation nodes, in which abnormal reductions in interpolation errors were obtained.

For 9 knots

$$1, 2, 3, 4, 5, 6, 7, 8, 9, 9.525, 7.28, 2.665, .46$$

For 13 knots

$$1, 2, 3, 4, 5, 6, 7, 8, 9, 10, 11, 12, 13, 13.37, 10.15, 3.846, .626$$

*Error tables for interpolation of function* $y = \sin(.75t)$, $x \in [0, 2\pi)$

Table 3. The ratio of interpolation errors without phantom nodes
to the errors of interpolation with 2 phantom nodes

| Number of Interpolation points | Continuity of the function | Continuity of the function and the first derivative | Continuity of the function and the second derivative | Values are selected |
|---|---|---|---|---|
| 5 | 0.73 | 1.5 | 4.8 | 10.6 |
| 9 | 5.8 | 50 | 84.6 | 220 |
| 13 | 4.4 | 23.1 | 60 | 133.3 |

Table 4. The ratio of interpolation errors without phantom nodes
to the errors of interpolation with 4 phantom nodes

| Number of Interpolation points | Continuity of the function | Continuity of the function and the first derivative | Continuity of the function and the second derivative | Values are selected |
|---|---|---|---|---|
| 5 | 0.58 | 1.4 | 10 | 597.82 |
| 9 | 3.14 | 14.6 | 53.7 | 1264 |
| 13 | 8 | 70.6 | 282.4 | 3428.6 |

*Error tables for interpolation of function* $y = .02 \exp t$, $x \in [0, 2\pi)$

Table 5. The ratio of interpolation errors without phantom nodes
to the errors of interpolation with 2 phantom nodes

| Number of Interpolation points | Continuity of the function | Continuity of the function and the first derivative | Continuity of the function and the second derivative | Values are selected |
|---|---|---|---|---|
| 5 | 1.5 | 2.75 | 5.5 | 16.9 |
| 9 | 1.8 | 13.75 | 15.7 | 23 |
| 13 | 2.1 | 5.2 | 12.1 | 27.1 |

Table 6. The ratio of interpolation errors without phantom nodes
to the errors of interpolation with 4 phantom nodes

| Number of Interpolation points | Continuity of the function | Continuity of the function and the first derivative | Continuity of the function and the second derivative | Values are selected |
|---|---|---|---|---|
| 5 | 1.6 | 5.2 | 8.9 | 49 |
| 9 | 2.2 | 7.9 | 22 | 407 |
| 13 | 2.1 | 8.1 | 8.9 | 506 |

Note that the value of the selected phantom points were looking by us with accuracy up to $10^{-3}$. By increasing this accuracy, one can get an increasing the coefficients shown in the tables.

Analyzing the data presented in the tables 1, 2, it is possible to put forward the hypothesis that at selected values of phantom points we obtain certain analogs of the trigonometric polynomials of the best approximation in space $C_{[0, 2\pi-\alpha]}$. However, it is difficult to prove this fact, using the well-known results of the approximation theory. In our opinion, this is due to the fact that in the theory of approximations, the number of points in the Chebyshev alternate determines the order of approximating polynomial; in our case, such a connection is violated.

If this hypothesis is correct, the task of finding points is greatly simplified, that are similar to points of "Chebyshevsky alternative". Indeed, in this case, the arrangement of these points at the segment is fixed; it is necessary to determine only values in a small number of phantom points.

In any case, the proposed method of phantom knots requires further research.

**Conclusions.**

There is proposed a method for improving the convergence of Fourier series by function systems, orthogonal at the segment, the application of which allows for smooth functions to receive uniformly convergent series. There is also proposed the method of phantom nodes improving the convergence of interpolation polynomials on systems of orthogonal functions, the application of which in many cases can significantly reduce the interpolation errors of these polynomials. The results of calculations are given at test cases using the proposed methods for trigonometric Fourier series; these calculations illustrate the high efficiency of these methods. Undoubtedly, the proposed method of phantom knots requires further theoretical studies.

**List of references**